\documentclass[11pt]{amsart}
\usepackage{amsmath,amssymb,amsthm}
\usepackage[a4paper, margin=2cm]{geometry}
\usepackage[backend=biber,style=alphabetic,giveninits=true,maxnames=99,minalphanames=3,sortcites=true,sorting=alphyear,date=year,doi=false,url=false,isbn=false]{biblatex}
\DeclareSortingTemplate{alphyear}{\sort{\field{presort}}\sort[final]{\field{sortkey}}\sort{\field{labelalpha}}\sort{\field{year}}\sort{\field{title}}}
\DeclareFieldFormat{postnote}{#1}
\DeclareFieldFormat{multipostnote}{#1}
\DeclareLabelalphaTemplate{\labelelement{\field[final]{shorthand}\field{label}\field[strwidth=3,strside=left,ifnames=1]{labelname}\field[strwidth=1,strside=left]{labelname}}}
\DeclareSourcemap{\maps[datatype=bibtex]{%
  \map[overwrite=true]{\step[fieldsource=shortjournal]\step[fieldset=journaltitle,origfieldval]}%
  \map{\pertype{online}\pertype{misc}\step[fieldset=date,null]\step[fieldset=year,null]}}}
\DeclareNameAlias{default}{given-family}

\DeclareFieldFormat*{title}{\mkbibemph{#1}}
\DeclareFieldFormat{journaltitle}{#1\isdot}
\DeclareFieldFormat[article]{volume}{\mkbibbold{#1}}
\DeclareFieldFormat{pages}{#1}
\DeclareFieldFormat{extraalpha}{#1}
\DeclareFieldFormat{pagetotal}{#1 pp\adddot}
\DeclareNolabel{\nolabel{\regexp{[\p{Pc}\p{Ps}\p{Pe}\p{Pi}\p{Pf}\p{Po}\p{S}\p{C}\p{M}]+}}}
\DeclareFieldFormat{eprint:arxiv}{arXiv:#1}
\AtEveryBibitem{\ifentrytype{incollection}{\clearfield{booktitle}}{}}
\renewbibmacro{in:}{\ifentrytype{article}{}{\ifentrytype{incollection}{}{\printtext{in}\setunit{\addspace}}}}
\renewbibmacro*{volume+number+eid}{\printfield{volume}}
\renewbibmacro*{issue+date}{\iffieldundef{year}{}{\printtext[parens]{\usebibmacro{date}}}\newunit}
\renewbibmacro*{note+pages}{\setunit{\addcomma\space}\printfield{note}\setunit{\addcomma\space}\printfield{eid}\setunit{\bibpagespunct}\printfield{pages}\setunit{\addcomma\space}\printfield{pagetotal}\newunit}
\renewbibmacro*{publisher+location+date}{\printlist{publisher}\setunit*{\addcomma\space}\printlist{location}\setunit*{\addcomma\space}\usebibmacro{date}\newunit}
\renewbibmacro*{series+number}{\printfield{series}\setunit*{\addspace}\printfield{number}\newunit}
\renewbibmacro*{doi+eprint+url}{\ifentrytype{online}{\usebibmacro{eprint}}{}}
\appto{\bibsetup}{\emergencystretch=2em}

\usepackage{hyperref}
\usepackage{xcolor}
\definecolor{CharlesRed}{HTML}{D22D40}
\hypersetup{
  pdfauthor={Vitezslav Kala, Scott Duke Kominers, Martí Oller, Pavlo Yatsyna},
  colorlinks=true,
  allcolors=CharlesRed
}
\newtheorem{theorem}{Theorem}[section]
\newtheorem{corollary}[theorem]{Corollary}
\newtheorem{lemma}[theorem]{Lemma}
\newtheorem{proposition}[theorem]{Proposition}
\newtheorem{conjecture}[theorem]{Conjecture}
\theoremstyle{definition}

\theoremstyle{remark}
\newtheorem{remark}[theorem]{Remark}
\numberwithin{equation}{section}

\makeatletter
\def\subsection{\@startsection{subsection}{2}%
	\z@{.5\linespacing\@plus.7\linespacing}{.3\linespacing}%
	{\normalfont\bfseries}}
\makeatother

\newcommand{\Q}{\mathbb{Q}}
\newcommand{\Z}{\mathbb{Z}}
\newcommand{\R}{\mathbb{R}}
\newcommand{\OO}{\mathcal{O}}
\newcommand{\Tr}{\operatorname{Tr}}
\newcommand{\N}{\operatorname{N}}
\newcommand{\cM}{\mathcal{M}}

\newcommand{\fp}{\mathfrak p}
\newcommand{\fP}{\mathfrak P}

\newcommand{\sgn}{\operatorname{sgn}}

\newcommand{\eps}{\varepsilon}

\makeatletter
\newcommand{\itemlab}[2]{\item[\textup{#1}]\protected@edef\@currentlabel{#1}\phantomsection\label{#2}}
\newcommand{\hypitem}[2]{\item[\textbf{#1}]\protected@edef\@currentlabel{#1}\phantomsection\label{#2}}
\makeatother
\newcommand{\partref}[2]{\hyperref[#1:#2]{\ref*{#1}\ref*{#1:#2}}}

\title{Kitaoka's Conjecture on universal quadratic forms is true}

\author{V\' \i t\v ezslav Kala}
\author{Scott Duke Kominers}
\author{Mart\'i Oller}
\author{Pavlo Yatsyna}

\address[Kala, Oller, and Yatsyna]{Charles University, Faculty of Mathematics and Physics, Department of Algebra, Sokolov\-sk\' a 83, 18600 Praha 8, Czech Republic}
\email[Kala]{vitezslav.kala@matfyz.cuni.cz}
\email[Oller]{marti.oller-riera@matfyz.cuni.cz}
\email[Yatsyna]{p.yatsyna@matfyz.cuni.cz}
\address[Kominers]{Harvard Business School; Department of Economics and Center of
Mathematical Sciences and Applications, Harvard University; and a16z Frontier Research}
\email[Kominers]{kominers@fas.harvard.edu}


\begin{document}

	\begin{abstract}
		We prove that the only totally real number fields admitting a universal classical ternary quadratic form are $\Q(\sqrt2)$,
		$\Q(\sqrt3)$, and $\Q(\sqrt5)$, settling a conjecture by Kitaoka from the early 1990s.
	\end{abstract}

	\maketitle

	\section{Introduction}

	In 1770, Lagrange proved that every positive integer is a sum of four squares, i.e., that the quadratic form $x^2+y^2+z^2+w^2$ is universal.
	On the other hand, no positive definite ternary quadratic form over $\Z$ is universal, already for local reasons.

	Over number fields, however, ternary quadratic forms can be universal.
	To be precise, let $K$ be a totally real number field with ring of integers $\OO_K$, and let $\OO_K^+$ be the set of totally positive elements of $\OO_K$, i.e., of the elements that are positive in every real embedding of $K$.
	A quadratic form $Q(x_1,\dots,x_r)=\sum_{1\le i\le j\le r}a_{ij}x_ix_j$ with $a_{ij}\in\OO_K$ is \emph{positive definite} if $Q(x)\in\OO_K^+$ for every nonzero $x\in\OO_K^r$, and \emph{classical} if $a_{ij}\in2\OO_K$ for all $i\ne j$, i.e., if its Gram matrix has entries in $\OO_K$.
	A positive definite quadratic form is \emph{universal} over $K$ if it represents every element of $\OO_K^+$.
	In 1941, Maa\ss{} \cite{Ma} proved that $x^2+y^2+z^2$ is universal over $\Q(\sqrt5)$.
	In the early 1990s, Kitaoka formulated the following conjecture.

	\begin{conjecture}[Kitaoka's Conjecture]
		There exist only finitely many totally real number fields admitting a universal classical ternary quadratic form.
	\end{conjecture}

	For more than thirty years, Kitaoka's Conjecture has been the main open problem on universal quadratic forms over number fields.
		In 1996, Chan--Kim--Raghavan \cite[Theorem 1.1]{CKR} proved that $\Q(\sqrt2)$, $\Q(\sqrt3)$, and $\Q(\sqrt5)$ are the only real quadratic fields with a universal classical ternary quadratic form, and they determined all such forms up to isometry, namely
	\begin{equation}\label{eq:CKR}
	\begin{aligned}
		\Q(\sqrt2)&\colon\ x^2+y^2+(2+\sqrt2)z^2,\ \ x^2+(2+\sqrt2)y^2+2yz+(2-\sqrt2)z^2,\\
		&\phantom{\colon}\ x^2+(2\pm\sqrt2)y^2+2yz+3z^2,\\
		\Q(\sqrt3)&\colon\ x^2+y^2+(2+\sqrt3)z^2,\ \ x^2+(2+\sqrt3)(y^2+z^2),\\
		\Q(\sqrt5)&\colon\ x^2+y^2+z^2,\ \ x^2+y^2+2z^2,\ \ x^2+y^2+\tfrac{5+\sqrt5}2z^2,\ \ x^2+2y^2+2yz+\tfrac{5\pm\sqrt5}2z^2.
	\end{aligned}
	\end{equation}
    
	Further partial results include the following:
	\begin{itemize}
		\item If the degree $[K:\Q]$ is odd, then no ternary quadratic form is universal over $K$ \cite[Lemma~3]{EK}.
		\item For every $d$, only finitely many totally real number fields of degree $d$ admit a universal ternary quadratic form \cite[Theorem 3]{KY23}.
		For classical forms, this was proven earlier by B. M. Kim in an unpublished manuscript.
		\item No totally real number field of degree $4$ admits a universal classical ternary form \cite[Theorem 1.1]{KK26}, which generalizes the case of biquadratic fields \cite{KTZ}.
		\item If $2$ is unramified in $K$ and $K$ admits a universal classical ternary form, then $K=\Q(\sqrt5)$ \cite[Theorem 1.1]{KKK}.
		\item The only real cyclotomic fields $\Q(\zeta_k+\zeta_k^{-1})$, where $\zeta_k=e^{2\pi i/k}$, with a universal classical ternary form are $\Q(\sqrt2)$, $\Q(\sqrt3)$, and $\Q(\sqrt5)$ \cite[Theorem 1]{Krasensky}, extending the case when $k$ is a prime or a power of $2$ \cite[Theorem 1.4]{KY24}.
	\end{itemize}

	We prove Kitaoka's Conjecture in the following explicit form.

	\begin{theorem}\label{thm:main}
		Let $K$ be a totally real number field.
		Some classical ternary quadratic form is universal over $K$ if and only if $K = \Q(\sqrt2)$, $\Q(\sqrt3)$, or $\Q(\sqrt5)$.
	\end{theorem}

	In particular, \eqref{eq:CKR} is the complete list of universal classical ternary forms over all totally real number fields.
	The ``if'' direction is due to Chan--Kim--Raghavan \cite{CKR}, and the ``only if'' direction was conjectured by Kramer--Kr\'asensk\'y \cite[Conjecture 1.2]{KK26}.
	Their conjecture also covers classical ternary lattices that are not free, i.e., that do not correspond to quadratic forms (see Section \ref{sec:prelim}).
	In this generality, we also prove the statement for fields in which some totally positive unit is not the square of a unit (Corollary \ref{cor:lattices}).
		
	\medskip

	Kitaoka's Conjecture is part of the broader study of universal quadratic forms and their minimal ranks.
	Over $\Z$, universal forms are completely characterized by the celebrated $15$- and $290$-Theorems of Conway--Schneeberger and Bhargava--Hanke \cite{Bh,BH}.
	Over every totally real number field, universal forms exist by the asymptotic local-global principle of Hsia--Kitaoka--Kneser \cite{HKK}, but their ranks can be large.
	Kala--Yatsyna--\.Zmija \cite{KYZ} showed that for every $r$, the real quadratic fields with a universal form of rank at most $r$ have density $0$ (see also \cite{BK,BK18,Ka16}), and some similar but weaker results hold in higher degrees \cite{Ka23,KT,Man,Yat19}.
	Among real quadratic fields, Kim--Kim--Park \cite{KKP} proved that only finitely many admit a universal form of rank at most $7$, while B. M. Kim \cite{Kim} showed that infinitely many admit one of rank $8$.
	Over fields of higher degree, the analogous questions for small ranks remain wide open.
	For related results and further background, see, e.g., \cite{COh,HHX,XZ} and the surveys \cite{KalaSurvey,Kim04}.

	\medskip

	Our approach is inspired by Siegel's proof \cite{Si-sq} that $\Q$ and $\Q(\sqrt5)$ are the only totally real number fields over which every element of $\OO_K^+$ is a sum of squares (see \cite[Theorem 8.1]{KalaSurvey} for an overview).
	In modern terms, his argument amounts to showing that these are the only fields in which all \emph{indecomposables} are squares, where an element $\alpha\in\OO_K^+$ is indecomposable if it cannot be written as $\alpha=\beta+\gamma$ with $\beta,\gamma\in\OO_K^+$.
    	To this end, he chose a suitable subfield $F\subsetneq K$ and considered the elements of $\OO_K^+\setminus F$ of minimal trace, which we call \emph{$F$-minimal}.
		These elements are indecomposable (Lemma \partref{lem:minimizer}{a}), and so they are squares if all indecomposables are.

	Our proof reduces Theorem \ref{thm:main} to statements about indecomposables, some of which, such as Theorem \ref{thm:HS} and Corollary \ref{cor:sqclass}, may be of independent interest.
	By \cite[Theorem 1.1]{KKK}, we may assume that $2$ is ramified in $K$.
	We will distinguish two cases, according to whether the group $U^+$ of totally positive units of $\OO_K$ equals the group $U^2$ of squares of units.

	If $U^+\ne U^2$, then \cite[Proposition 2.5 and Theorem 3.1(d)]{KKK} show that a universal classical ternary form can exist only when every indecomposable is a square or a half-square, i.e., of the form $x^2$ or $x^2/2$ with $x\in\OO_K$.
	In Section \ref{sec:HS}, we will prove that this property forces $K=\Q$, $\Q(\sqrt3)$, or $\Q(\sqrt5)$ (Theorem \ref{thm:HS}).

	If $U^+=U^2$, let $Q$ be a universal classical ternary quadratic form over $K$.
	As $Q$ represents $1$, we have $Q\cong\langle1\rangle\perp Q_0$ for a binary form $Q_0$.
	In Section \ref{sec:lattices}, we will use a representation of the determinant of $Q_0$ by $Q$ to write $Q_0=\ell^2+R$, where $\ell$ is an integral linear form and $R$ is a positive semidefinite binary form whose determinant is a square.
	Then $R$ represents every nonsquare indecomposable $\alpha$, and $\alpha R(v)$ is a sum of two integral squares for every $v\in\OO_K^2$.
	Using these representations and a representation of $3$ by $Q$, we will establish two properties: all nonunit indecomposables lie in a single dyadic prime, and twice every indecomposable is a sum of two integral squares.
	In Section \ref{sec:dyadic}, we will prove that among fields with $U^+=U^2$ in which $2$ is ramified, only $\Q(\sqrt2)$ has these properties (Proposition \ref{prop:dyadic}), and we will also deduce that $\Q(\sqrt2)$ is the only field containing $\sqrt2$ in which every indecomposable has the form $x^2$ or $(2+\sqrt2)x^2$ with $x\in\OO_K$ (Corollary \ref{cor:sqclass}).
	
	\medskip

	The analogue of Kitaoka's Conjecture for nonclassical quadratic forms, i.e., without the condition $a_{ij}\in2\OO_K$, remains open, although only finitely many real quadratic fields admit such a universal ternary form \cite{KKP}, in fact at most $13$ \cite{KKPYZ}.
    Kr\'asensk\'y--Scharlau announced nonclassical universal ternary forms over the quartic fields $\Q(\sqrt2,\sqrt5)$ and $\Q(\zeta_{20}+\zeta_{20}^{-1})$, and so Theorem \ref{thm:main} does not extend to nonclassical forms.

	Hotta \cite{Hotta-I,Hotta-II} has announced a proof of Theorem \ref{thm:main} by a different method, assuming the generalized Riemann hypothesis for Dedekind zeta functions, whereas our proof is much shorter and unconditional.

	\subsection*{Acknowledgments}

    We thank the members of the \href{https://www.karlin.mff.cuni.cz/~kala/web/ufoclan}{UFOCLAN} research group for many enjoyable and helpful discussions over the past ten years, especially Jakub Kr\'asensk\'y and Robin Visser for comments on drafts of this manuscript.
	This work was conducted while S.D.K. was visiting the Technological Innovation, Entrepreneurship, and Strategic Management
(TIES) Group at the MIT Sloan School of Management; he greatly appreciates their hospitality.
	V.K., M.O., and P.Y. were supported by Czech Science Foundation (GA\v{C}R) grant 26-20514S, and P.Y. by Charles University programme PRIMUS/24/SCI/010.

	\subsection*{Use of generative AI}

	An initial proof and draft were produced with the assistance of Claude (Anthropic), building on unpublished ideas of the authors and the works cited here.
	Further refinements of the argument were explored with GPT-6.2 Astra (OpenAI).
	We thank Mira Dot (OpenAI) for helpful research assistance, particularly in optimizing the determinant argument, and for many conversations about alternative readings of the proof.

	\section{Preliminaries}\label{sec:prelim}

	Throughout this paper, $K$ is a totally real number field of degree $d$ with embeddings $\sigma_1,\dots,\sigma_d\colon K\hookrightarrow\R$, so that the conjugates of $\alpha\in K$ are $\sigma_1(\alpha),\dots,\sigma_d(\alpha)$.
	We write $\N$ and $\Tr$ for the norm and trace from $K$ to $\Q$, i.e., $\N(\alpha)=\prod_i\sigma_i(\alpha)$ and $\Tr(\alpha)=\sum_i\sigma_i(\alpha)$.
	For $\alpha,\beta\in K$, we say that $\alpha$ is \emph{totally greater than} $\beta$, and write $\alpha\succ\beta$ or $\beta\prec\alpha$, if $\sigma_i(\alpha)>\sigma_i(\beta)$ for all $i$.
	We write $\alpha\succeq\beta$ or $\beta\preceq\alpha$ if $\alpha\succ\beta$ or $\alpha=\beta$.
	An element $\alpha\in K$ is \emph{totally positive} if $\alpha\succ0$, and the \emph{signature} of $\alpha\in K^\times$ is $\sgn(\alpha)=(\sgn\sigma_1(\alpha),\dots,\sgn\sigma_d(\alpha))\in\{\pm1\}^d$.

	Let $\OO=\OO_K$ be the ring of integers of $K$, and let $\OO^+$ be the set of totally positive elements of $\OO$.
	We denote by $U=\OO^\times$ the group of units, by $U^+=U\cap\OO^+$ its subgroup of totally positive units, and by $U^2=\{u^2:u\in U\}$ its subgroup of squares.
		For a totally real number field $F$, we similarly write $\OO_F$, $U_F$, $U_F^+$, and $U_F^2$ for the corresponding objects of $F$, while the notation without a subscript always refers to $K$.

	For a nonzero prime ideal $\fP$ of $\OO$, we denote by $v_\fP$ the $\fP$-adic valuation on $K$, normalized so that $v_\fP(K^\times)=\Z$.
	If $\fP$ lies above $2$, we call $\fP$ a \emph{dyadic} prime and denote by $e_\fP=v_\fP(2)$ its ramification index.
	We put $\varphi=(1+\sqrt5)/2$, so that $\varphi^{\pm2}=(3\pm\sqrt5)/2$.

	An element of $\OO$ is a \emph{square} if it equals $x^2$ with $x\in\OO$, and a \emph{half-square} if it equals $x^2/2$ with $x\in\OO$.
	Since $\OO$ is integrally closed,  the condition $x \in \OO$ can be relaxed to $x \in K$.

	\medskip

	A \emph{quadratic form} over $\OO$ of rank $r$ is a homogeneous polynomial $Q(x_1,\dots,x_r)=\sum_{1\le i\le j\le r}a_{ij}x_ix_j$ with $a_{ij}\in\OO$.
	Its \emph{Gram matrix} is the symmetric matrix $G$ with diagonal entries $a_{ii}$ and off-diagonal entries $a_{ij}/2$, so that $Q(x)=x^{\mathsf T}Gx$.

	A \emph{quadratic lattice} over $\OO$ is a pair $(L,Q)$ consisting of a finitely generated torsion-free $\OO$-module $L$ and a map $Q \colon L \to \OO$ such that
	\begin{enumerate}
		\item[(1)] $Q(\alpha v) = \alpha^2 Q(v)$ for all $\alpha \in \OO$ and $v \in L$, and
		\item[(2)] the map $B \colon L \times L \to K$ defined by $B(u,v) = \frac{1}{2}(Q(u+v)-Q(u)-Q(v))$ is bilinear.
	\end{enumerate}	
    We will call a quadratic lattice $(L,Q)$ simply a \emph{lattice} and often denote it just by $L$.
	The \emph{rank} of $L$ is its rank as an $\OO$-module, and $L$ is \emph{binary} if its rank is $2$ and \emph{ternary} if its rank is $3$.    
    Furthermore, a lattice $L$ is \emph{positive definite} if $Q(v) \in \OO^+$ for all nonzero $v \in L$, \emph{positive semidefinite} if $Q(v)\in\OO^+\cup\{0\}$ for all $v\in L$, and \emph{classical} if $B(u,v) \in \OO$ for all $u,v \in L$.
	A lattice $L$ \emph{represents} $\alpha\in\OO$ if $Q(v)=\alpha$ for some $v\in L$, and $L$ is \emph{universal} if it represents every element of $\OO^+$.	We extend $Q$ and $B$ to the $K$-vector space $L\otimes K=L\otimes_\OO K$, which we call the \emph{quadratic space} of $L$.

	Two lattices $(L,Q)$ and $(L',Q')$ are \emph{isometric}, written $L\cong L'$, if there is an $\OO$-linear bijection $\iota\colon L\to L'$ with $Q'(\iota(v))=Q(v)$ for all $v\in L$.
	The \emph{orthogonal sum} $L\perp L'$ is the $\OO$-module $L\oplus L'$ with the quadratic map $(v,v')\mapsto Q(v)+Q'(v')$.

	A lattice is \emph{free} if it is a free $\OO$-module.
	A quadratic form $Q$ of rank $r$ corresponds to the free lattice $(\OO^r,Q)$, and this gives a bijection between quadratic forms up to invertible linear changes of variables over $\OO$ and free lattices up to isometry.
	Accordingly, we apply the notions defined for lattices also to quadratic forms.
	For example, $Q$ is classical if and only if its Gram matrix has entries in $\OO$, i.e., $a_{ij}\in2\OO$ for all $i\ne j$.
	Given $\alpha_1,\dots,\alpha_n\in\OO^+$, $\langle\alpha_1,\dots,\alpha_n\rangle$ denotes the diagonal form $\sum\alpha_ix_i^2$.
	Unless stated otherwise, all quadratic forms over $\OO$ and all lattices in this paper are classical and positive definite.

	\medskip

	An element $\alpha \in \OO^+$ is \emph{indecomposable} if it cannot be written as $\alpha = \beta + \gamma$ with $\beta,\gamma \in \OO^+$; it is \emph{decomposable} otherwise.

	\begin{lemma}\label{lem:basic}
		Let $\alpha,\beta\in\OO^+$.
		Then $\N(\alpha+\beta)^{1/d}\ge\N(\alpha)^{1/d}+\N(\beta)^{1/d}$, and $\Tr(\alpha)\ge d\,\N(\alpha)^{1/d}\ge d$, with $\Tr(\alpha)=d$ only for $\alpha=1$.
		In particular, every $\alpha\in\OO^+$ with $\N(\alpha)<2^d$ is indecomposable.
	\end{lemma}

	\begin{proof}
		The inequality for the norm is well known (see, e.g., \cite[3.1]{OM80}), and it shows that every decomposable element of $\OO^+$ has norm at least $2^d$.
		The inequality for the trace follows from the AM-GM inequality and the fact that  $\N(\alpha)\ge1$.
	\end{proof}

	The group $\sgn(U)$ of signatures of units has index $|U^+/U^2|$ in $\{\pm1\}^d$ (see, e.g., \cite[Lemma 3.10(a)]{KKK}).
	In particular, if $U^+=U^2$, then there are units of all signatures.

	\section{Small elements and indecomposables}\label{sec:small}

	In this section we collect a number of results on indecomposables that will be useful in our later arguments; these results do not involve quadratic forms directly.

	The following lemma follows from a basic result on the Schur--Siegel--Smyth trace problem \cite[Theorem III]{Si-tr} (see also, e.g., \cite[Lemma 2.2]{KKK}).

	\begin{lemma}\label{lem:small}
		If $w\in\OO$ and $0\preceq w\preceq3$, then $w\in\{0,1,2,3,\varphi^2,\varphi^{-2}\}$.
	\end{lemma}

	We next give a criterion for deducing the indecomposability of $Y$ from that of $Y^2/A$ when $A$ is a totally positive divisor of $2$.
	
	\begin{lemma}\label{lem:SU}
		Let $Y\in\OO^+$ be such that $Y^2/A$ is an indecomposable element of $\OO$ for some $A\in\OO^+$ with $2/A\in\OO$.
		Then $Y$ is indecomposable.
	\end{lemma}

	\begin{proof}
		If $Y=\beta+\gamma$ with $\beta,\gamma\in\OO^+$, then
		\[
		\frac{Y^2}{A}=\frac{\beta^2+\gamma^2}{A}+\frac2A\,\beta\gamma .
		\]
		The second summand is integral because $2/A$ is, the first is integral as the difference of two elements of $\OO$, and both are totally positive.
		This contradicts the indecomposability of $Y^2/A$.
	\end{proof}
    
	We next show that a prime ideal containing all nonsquare indecomposables contains all nonunit indecomposables.

	\begin{lemma}\label{lem:notinP}
		Assume that $U^+=U^2$, and assume there exists a prime ideal $\fP$ of $\OO$ that contains every indecomposable that is not a square.
		Then every indecomposable that does not lie in $\fP$ is a unit.
	\end{lemma}

	\begin{proof}
		Let $\alpha\notin\fP$ be a nonunit indecomposable.
		As $U^+=U^2$, \cite[Lemma 3.5]{KK26} gives $\alpha=\mu\beta^{2^k}$ with $\mu\in U^+$, $k\ge0$, and an indecomposable $\beta$ that is not a square.
		By our assumption, $\beta\in\fP$, and hence $\alpha\in\fP$, a contradiction.
	\end{proof}

	\medskip

	As in Siegel's proof \cite{Si-sq}, we will consider elements of minimal trace outside a subfield of $K$, and more generally outside an additive subgroup of $K$.
	For an additive subgroup $S$ of $K$ with $1\in S$ and $\OO\not\subset S$, we put
	\[
	q_S=\min\{\Tr(x):x\in\OO^+\setminus S\}.
	\]
	Here $\OO^+\setminus S$ is nonempty, as $\beta+m\in\OO^+\setminus S$ for $\beta\in\OO\setminus S$ and all sufficiently large $m\in\Z$.
	An element $x\in\OO^+\setminus S$ is \emph{$S$-minimal} if $\Tr(x)=q_S$.

	\begin{lemma}\label{lem:minimizer}
		Let $S$ be an additive subgroup of $K$ with $1\in S$ and $\OO\not\subset S$. Then:
		\begin{enumerate}
			\itemlab{(a)}{lem:minimizer:a} Every $S$-minimal element is indecomposable.
			\itemlab{(b)}{lem:minimizer:b} Every decomposable $x\in\OO^+\setminus S$ satisfies $\Tr(x)\ge q_S+d$, and if $\Tr(x)=q_S+d$, then $x-1\in\OO^+$.
		\end{enumerate}
	\end{lemma}

	\begin{proof}
		(a) In a decomposition $x=\beta+\gamma$ of an $S$-minimal element $x$, one summand lies outside $S$, as $S$ is an additive group, and it has smaller trace.

		(b) As $S$ is an additive group, we can write $x=y+z$ with $y,z\in\OO^+$ and $y\notin S$.
		Then $\Tr(x)=\Tr(y)+\Tr(z)\ge q_S+d$ by Lemma \ref{lem:basic}, and equality forces $\Tr(z)=d$, so $z=1$ and $x-1=y\in\OO^+$.
	\end{proof}

	Next, we will consider $S$-minimal elements that are squares or half-squares.

	\begin{lemma}\label{lem:Z}
		Let $S$ be an additive subgroup of $K$ with $1\in S$ and $\OO\not\subset S$, and let $y\in\OO$ and $A\in\{1,2\}$ be such that $x=y^2/A$ is $S$-minimal.
		Then the elements $z_\pm=x\pm y+1$ are totally positive and are not units, and at least one of them lies in $S$ or is indecomposable.
	\end{lemma}

	\begin{proof}
		As $x\notin S$ and $\Z\subset S$, we have $x\notin\Z$, so $y\notin\Z$, and hence no conjugate of $y$ is rational.
		For $A=1$, every conjugate of $z_\pm=(y\pm\frac12)^2+\frac34$ is at least $3/4$, and for $A=2$, every conjugate of $z_\pm=((y\pm1)^2+1)/2$ is greater than $1/2$.
		So $|\sigma(z_\pm-1)|<\sigma(z_\pm)$ for every embedding $\sigma$, and hence $|\N(z_\pm-1)|<\N(z_\pm)$.
		If $z_+$ or $z_-$ were a unit, it would therefore be $1$, so $y(y\pm A)=0$ and $y\in\Z$, a contradiction.

		Suppose that $z_+,z_-\notin S$ and that both are decomposable.
		Then $\Tr(z_\pm)\ge q_S+d$ by Lemma \partref{lem:minimizer}{b}, while $\Tr(z_+)+\Tr(z_-)=2\Tr(x)+2d=2(q_S+d)$.
		So $\Tr(z_\pm)=q_S+d$, and Lemma \partref{lem:minimizer}{b} gives $x\pm y=z_\pm-1\succ0$.
		Multiplying, we get $x(x-A)=x^2-y^2\succ0$, so $x\succ A\ge1$, and $x=(x-1)+1$ is decomposable, contrary to Lemma \partref{lem:minimizer}{a}.
	\end{proof}

	For a dyadic prime $\fP$ with $e_\fP\ge2$, we have $2\in\fP^2$, so squaring is an additive map on $\OO/\fP^2$, and the set
	\[
	S_\fP=\{a\in\OO: a\equiv w^2\pmod{\fP^2}\text{ for some }w\in\OO\}
	\]
	is an additive subgroup of $\OO$ that contains $1$ and $\fP^2$.
	Moreover,
	\begin{equation}\label{eq:squareres}
		S_\fP\cap\fP=\fP^2,
	\end{equation}
	since if $a\in\fP$ and $a\equiv w^2\pmod{\fP^2}$ with $w\in\OO$, then $w\in\fP$, and so $a\in\fP^2$.
	In particular, $S_\fP\ne\OO$.
	Also, if $r,s\in\OO$ satisfy $r^2+s^2\in\fP$, then $r^2+s^2\in\fP^2$ by \eqref{eq:squareres}, as $r^2+s^2\equiv(r+s)^2\pmod{\fP^2}$.

	\begin{lemma}\label{lem:localtranslate}
		Let $\fP$ be a dyadic prime with $e_\fP\ge2$, and assume that there is an additive subgroup $G$ of $\OO/\fP^2$ with $1\notin G$ that contains the residue modulo $\fP^2$ of every indecomposable that is not a square.
		Then no $S_\fP$-minimal element has the form $y^2/2$ with $y\in\fP^2$.
	\end{lemma}

	\begin{proof}
		Suppose that $x=y^2/2$ is $S_\fP$-minimal with $y\in\fP^2$.
		As $x\notin S_\fP$ and $1,y\in S_\fP$, the elements $x\pm y+1$ do not lie in $S_\fP$.
		So one of them, say $z$, is indecomposable by Lemma \ref{lem:Z} with $A=2$.
		As $x,z\notin S_\fP$, neither $x$ nor $z$ is a square, and $x$ is indecomposable by Lemma \partref{lem:minimizer}{a}.
		So the residues of $x$ and $z$ modulo $\fP^2$ lie in $G$, and hence so does the residue of $z-x=1\pm y$, which is $1$, as $y\in\fP^2$.
		This contradicts $1\notin G$.
	\end{proof}

	When there are units of all signatures, we can also find indecomposables by minimizing the norm instead of the trace.

	\begin{lemma}\label{lem:leastnorm}
		Assume that $U^+=U^2$. Then:
		\begin{enumerate}
			\itemlab{(a)}{lem:leastnorm:a} If $\cM\subset\OO$ is a nonempty set with $\cM\cap2\OO=\emptyset$, $\cM+2\OO\subset\cM$, and $U\cM\subset\cM$, then $\cM\cap\OO^+$ is nonempty, and every element of $\cM\cap\OO^+$ of minimal norm is indecomposable.
			\itemlab{(b)}{lem:leastnorm:b} If $\fP$ is a dyadic prime with $e_\fP\ge2$, then there is an indecomposable $\gamma\in\OO^+$ with $v_\fP(\gamma)=1$ such that $\N(\gamma)\le\N(\gamma')$ for all $\gamma'\in\OO^+$ with $v_\fP(\gamma')=1$.
		\end{enumerate}
	\end{lemma}

	\begin{proof}
		(a) For $\beta\in\cM$ and all sufficiently large $m\in\Z$, we have $\beta+2m\in\cM\cap\OO^+$.
		Let $\gamma\in\cM\cap\OO^+$ be of minimal norm, and suppose that $\gamma=\beta_1+\beta_2$ with $\beta_1,\beta_2\in\OO^+$.
		Then $\beta_1-\beta_2=\gamma-2\beta_2\in\cM$, so $\beta_1-\beta_2\ne0$.
		As $|\sigma(\beta_1-\beta_2)|<\sigma(\gamma)$ for every embedding $\sigma$, we get $|\N(\beta_1-\beta_2)|<\N(\gamma)$.
		As $U^+=U^2$, there is $\eta\in U$ with $\eta(\beta_1-\beta_2)\succ0$.
		Then $\eta(\beta_1-\beta_2)\in\cM\cap\OO^+$ has smaller norm than $\gamma$, a contradiction.

		(b) The set $\cM=\fP\setminus\fP^2$ of elements of $\OO$ of valuation $1$ at $\fP$ is nonempty, and it satisfies the assumptions of \ref{lem:leastnorm:a}, as $2\OO\subset\fP^2$.
	\end{proof}

\section{Universal forms and sums of two squares}\label{sec:lattices}

	In this section we will derive properties of $K$'s indecomposables from the hypothesized existence of a universal classical ternary form.
	If there is such a form $Q$, we denote by $L$ the corresponding free lattice over $\OO$.

	As $L$ represents $1$, there must be $u_0\in L$ with $Q(u_0)=1$, and as $L$ is classical, we have $L=\OO u_0\perp M$ for the binary lattice $M=\{w\in L:B(u_0,w)=0\}$ (see, e.g., \cite[82:15]{OM}).
	Since the Steinitz class is additive and $L$ is free, $M$ must be free as well.
	Given this freeness, a global basis of $M$ will yield an $\OO$-valued alternating form $[\,\cdot\,,\,\cdot\,]$ and a totally positive determinant $\Delta\in\OO$ to which we can apply the universality hypothesis.
	
	We denote by $Q_0$ and $B$ the quadratic map and the bilinear form of $M$, and also their extensions to the quadratic space $V=M\otimes K$.
	We fix a basis $b_1,b_2$ of the $\OO$-module $M$, and for vectors $v=v_1b_1+v_2b_2$ and $w=w_1b_1+w_2b_2$ in $V$, we put $[v,w]=v_1w_2-v_2w_1$.
	For a quadratic form $H$ on $V$ with bilinear form $B_H$, we put $\det H=\det\bigl(B_H(b_i,b_j)\bigr)_{i,j}$, and then
	\begin{equation}\label{eq:gram}
		H(v)H(w)-B_H(v,w)^2=(\det H)[v,w]^2\qquad\text{for all }v,w\in V.
	\end{equation}

	Now, the determinant $\Delta=\det Q_0$ lies in $\OO^+$. As $L$ is universal, we thus have $c\in\OO$ and $f\in M$ such that
	\[
	\Delta=c^2+Q_0(f).
	\]
	We fix those values and put
	\[
	\ell(v)=[v,f],\qquad R=Q_0-\ell^2,\qquad\omega(v,w)=c\,[v,w].
	\]
	In the rest of the section, we will use the following hypothesis:
	\begin{itemize}
		\hypitem{(L)}{hyp:L} There is a free universal classical ternary lattice $L$ over $\OO$, and $M$, $Q_0$, $B$, $V$, $[\,\cdot\,,\,\cdot\,]$, $\Delta$, $c$, $f$, $\ell$, $R$, and $\omega$ are as above.
	\end{itemize}
	
\begin{lemma}\label{lem:detres}
  Assume \textup{\ref{hyp:L}}. Then the form $R$ is classical and positive semidefinite, and
  $\det R=c^2$. More precisely, we have
  \begin{equation}\label{eq:res}
    \Delta R(v)=c^2Q_0(v)+B(v,f)^2 \qquad \text{for all $v\in V$},
  \end{equation}
  and also
  \begin{equation}\label{eq:Romega}
    R(v)R(w)=B_R(v,w)^2+\omega(v,w)^2 \qquad \text{for all $v,w\in V$}.
  \end{equation}
  If $c\ne0$, then $R$ is positive definite.
\end{lemma}

\begin{proof}
  Since $\ell(M)\subset\OO$ and $B_R(v,w)=B(v,w)-\ell(v)\ell(w)$, the form $R$ is classical.
  Applying \eqref{eq:gram} with $H=Q_0$ and $w=f$ gives
  \[
  \Delta\ell(v)^2=Q_0(v)Q_0(f)-B(v,f)^2;
  \]
  subtracting this identity from $\Delta Q_0(v)$ gives \eqref{eq:res}.
  In particular, $R(v)\succeq0$ for every $v\in V$, with $R(v)\succ0$ for $v\ne0$ if $c\ne0$.
	
  It remains to compute the determinant. If $f=0$, then $R=Q_0$ and $\det R=\Delta=c^2$.
	Otherwise, we choose $v\in V$ with $[f,v]\ne0$.
  Since $\ell(f)=0$, applying \eqref{eq:gram} first to $R$ and then to $Q_0$ gives
  \begin{align*}
    (\det R)[f,v]^2
      &=Q_0(f)R(v)-B(f,v)^2\\
      &=\Delta[f,v]^2-Q_0(f)[v,f]^2\\
      &=c^2[f,v]^2.
  \end{align*}
  Thus $\det R=c^2$, and \eqref{eq:Romega} follows from \eqref{eq:gram} for $R$.
\end{proof}

	The identity \eqref{eq:Romega} in Lemma \ref{lem:detres} shows that multiplying two integral values of $R$ gives a sum of two integral squares.
	We next show that one of these values may be any prescribed nonsquare indecomposable.

	\begin{lemma}\label{lem:res}
		Assume \textup{\ref{hyp:L}}, and let $\alpha$ be an indecomposable that is not a square.
		Then there is $e\in M$ with $\ell(e)=0$ and $Q_0(e)=R(e)=\alpha$, and every such $e$ satisfies
		\begin{equation}\label{eq:resid1}
			\alpha R(v)=B(e,v)^2+\omega(e,v)^2\qquad\text{for all }v\in V.
		\end{equation}
	\end{lemma}

	\begin{proof}
		As $L$ is universal, there are $X\in\OO$ and $e\in M$ with $\alpha=X^2+Q_0(e)=X^2+\ell(e)^2+R(e)$.
		The three summands lie in $\OO$, and each of them is zero or totally positive.
		As $\alpha$ is indecomposable and not a square, we get $X=\ell(e)=0$, and so $Q_0(e)=R(e)=\alpha$.
		For every such $e$, we have $B_R(e,v)=B(e,v)$ for $v\in V$, and so \eqref{eq:Romega} for the vectors $e$ and $v$ gives \eqref{eq:resid1}.
	\end{proof}

	We will first use Lemma \ref{lem:res} for an indecomposable of valuation $1$ at a dyadic prime.

	\begin{proposition}\label{prop:resP}
		Assume \textup{\ref{hyp:L}} and that $U^+=U^2$, and let $\fP$ be a dyadic prime with $e_\fP\ge2$. Then:
		\begin{enumerate}
			\itemlab{(a)}{prop:resP:a} We have $R(v)\in\fP$ for all $v\in M$, and every indecomposable that is not a square lies in $\fP$.
			\itemlab{(b)}{prop:resP:b} Every indecomposable that does not lie in $\fP$ is a unit, and $\fP$ is the only dyadic prime of $K$.
		\end{enumerate}
	\end{proposition}

	\begin{proof}
		(a) Let $\gamma$ be as in Lemma \partref{lem:leastnorm}{b}.
		As $v_\fP(\gamma)=1$, the element $\gamma$ is not a square, so Lemma \ref{lem:res} gives $e'\in M$ such that $\gamma R(v)=B(e',v)^2+\omega(e',v)^2$ for $v\in M$, where $B(e',v),\omega(e',v)\in\OO$.
		As $\gamma\in\fP$, this sum of two squares lies in $\fP$, and hence in $\fP^2$ (see \eqref{eq:squareres}).
		Since $v_\fP(\gamma)=1$, it follows that $R(v)\in\fP$.
		Finally, if $\alpha$ is an indecomposable that is not a square, then $\alpha=R(e)\in\fP$ for $e$ as in Lemma \ref{lem:res}.

		(b) By \ref{prop:resP:a} and Lemma \ref{lem:notinP}, every indecomposable that does not lie in $\fP$ is a unit.
		Suppose that $K$ has a dyadic prime $\fP'\ne\fP$.
		Then the set $\cM=\fP'\setminus\fP$ is nonempty, and it satisfies the assumptions of Lemma \partref{lem:leastnorm}{a}, as $2\OO\subset\fP\cap\fP'$.
		So Lemma \partref{lem:leastnorm}{a} gives a nonunit indecomposable in $\cM$, which does not lie in $\fP$, a contradiction.
	\end{proof}

	\begin{proposition}\label{prop:three}
		Assume \textup{\ref{hyp:L}}, that $U^+=U^2$, and that $2$ is ramified in $K$.
		Then $2\alpha$ is a sum of two squares of elements of $\OO$ for every indecomposable $\alpha$.
	\end{proposition}

	\begin{proof}
		If an indecomposable $\alpha$ is a square, then $2\alpha=\alpha+\alpha$ is a sum of two squares.
		So let $\alpha$ be an indecomposable that is not a square, and let $\fP$ be a dyadic prime with $e_\fP\ge2$.
		By Lemma \ref{lem:res}, choose $e\in M$ with $\ell(e)=0$ and $R(e)=\alpha$.
		We will write $2\alpha$ as a sum of two squares by representing $3$ by $L$, and in one case also $3a$ for a suitable $a\in\OO^+$.
		The element $3$ is not a square in $K$.
		Otherwise $K$ contains $F=\Q(\sqrt3)$, whose totally positive unit $2+\sqrt3$ is not a square in $F$, and \cite[Theorem 2.1]{EMP} gives $|U^+/U^2|\ge|U_F^+/U_F^2|\ge2$, contrary to $U^+=U^2$.

		As $L$ is universal, there are $X\in\OO$ and $g\in M$ with $3=X^2+Q_0(g)=X^2+\ell(g)^2+R(g)$, and $g\ne0$, as $3$ is not a square.
		As $0\preceq R(g)\preceq3$, Lemma \ref{lem:small} gives $R(g)\in\{0,1,2,3,\varphi^2,\varphi^{-2}\}$.
		By Proposition \partref{prop:resP}{a}, we also have $R(g)\in\fP$, and so $R(g)\in\{0,2\}$.

		\emph{Case $R(g)=2$.}
		Taking $v=g$ in \eqref{eq:resid1}, we get $2\alpha=B(e,g)^2+\omega(e,g)^2$, with both $B(e,g)$ and $\omega(e,g)$ in $\OO$.

		\emph{Case $R(g)=0$.}
		We will first show that there are $a\in\OO^+$ with $v_\fP(a)=1$ and an $\OO$-linear map $t\colon M\to\OO$ such that $R(v)=a\,t(v)^2$ for all $v\in M$.
		As $g\ne0$ and $R(g)=0$, we have $c=0$ by Lemma \ref{lem:detres}.
		So $Q_0(f)=\Delta$ and $\omega=0$, and \eqref{eq:res} gives $R(v)=B(v,f)^2/\Delta$ for $v\in V$.
		Let $\gamma$ be as in Lemma \partref{lem:leastnorm}{b}.
		As $v_\fP(\gamma)=1$, the element $\gamma$ is not a square.
		By Lemma \ref{lem:res}, choose $e'\in M$ with $\ell(e')=0$ and $R(e')=\gamma$.
		Since $\omega=0$, equation \eqref{eq:resid1} gives $\gamma R(v)=B(e',v)^2$ for all $v\in V$.

		If $\ell(g)=[g,f]$ is a unit, then $g$ and $f$ form a basis of $M$.
		As $R(g)=0$, we have $B(g,f)=0$, so $R(Yg+Zf)=Z^2R(f)$ for $Y,Z\in\OO$, and we take $a=R(f)=Q_0(f)=\Delta$ and $t(Yg+Zf)=Z$.
		Then $\gamma=R(e')=a\,t(e')^2$, so $v_\fP(a)=1-2v_\fP(t(e'))$ is odd and at most $1$, and hence $v_\fP(a)=1$.

		If $\ell(g)$ is not a unit, then $3=X^2+\ell(g)^2$ and $0\preceq\ell(g)^2\preceq3$, so Lemma \ref{lem:small} gives $\ell(g)^2\in\{0,2,3\}$.
		As $3$ is not a square, $X$ and $\ell(g)$ are nonzero, so $\ell(g)^2\prec3$, and hence $\ell(g)^2=2$.
		So $\sqrt2=\pm\ell(g)\in K$.

		As $U^+=U^2$, some associate of $\sqrt2$ is totally positive, and so we can choose a nonunit $\pi\in\OO^+$ that divides $\sqrt2$ such that $\N(\pi)$ is minimal among all such elements.
		By the minimality of $\N(\pi)$, the element $\pi$ is irreducible, and hence it is not a square.
		As $\N(\pi)\le|\N(\sqrt2)|=2^{d/2}<2^d$, the element $\pi$ is indecomposable by Lemma \ref{lem:basic}.
		By Lemma \ref{lem:res}, $\pi$ and every indecomposable $\beta$ that is not a square are values of $R$, so $\pi\Delta$ and $\beta\Delta$ are squares in $K$, and hence so is $\pi\beta$.
		As $\pi\beta\in\OO$, we get $\pi\beta=y^2$ with $y\in\OO$.
		Hence every indecomposable $\beta$ has the form $\beta=y^2/D$ with $D\in\{1,\pi\}$ and $y\in\OO$.

		Let us show that every indecomposable $\beta$ satisfies $v_\fp(\beta)=0$ for every prime ideal $\fp$ of $\OO$ with $2\notin\fp$.
		Write $\beta=y^2/D$ as above, and choose $\eta\in U$ with $\eta y\succ0$.
		As $(\eta y)^2/D=\eta^2\beta$ is indecomposable and $2/D\in\OO$, the element $\eta y$ is indecomposable by Lemma \ref{lem:SU} with $A=D$, and $v_\fp(\eta y)=v_\fp(\beta)/2$.
		Repeating this, we see that $v_\fp(\beta)$ is divisible by every power of $2$, and so $v_\fp(\beta)=0$.

		As $\fP$ is the only dyadic prime by Proposition \partref{prop:resP}{b}, we must have $(\gamma)=\fP$.
		Thus we may take $a=\gamma$ and $t(v)=B(e',v)/\gamma$.
		Indeed, $v_\fP(a)=1$, and Proposition \partref{prop:resP}{a} gives $R(v)\in\fP=(\gamma)$ for every $v\in M$.
		Thus $t(v)^2=R(v)/\gamma\in\OO$, so $t(v)\in\OO$.

		\medskip

		In either subcase, we now have $R=a\,t^2$ with $v_\fP(a)=1$.
		As $L$ is universal, there are $r\in\OO$ and $w\in M$ with $3a=r^2+Q_0(w)=r^2+\ell(w)^2+a\,t(w)^2$.
		Putting $s=\ell(w)$ and $z=t(w)$, we get
		\[
		a(3-z^2)=r^2+s^2\qquad\text{with }r,s,z\in\OO,
		\]
		and $0\preceq z^2\preceq3$.
		As $a\in\fP$, we have $r^2+s^2\in\fP$, so $r^2+s^2\in\fP^2$ (see \eqref{eq:squareres}).
		Since $v_\fP(a)=1$, this gives $z^2\equiv3\equiv1\pmod\fP$.
		
		Among the possible values $0,1,2,3,\varphi^2,\varphi^{-2}$ of $z^2$ given by Lemma \ref{lem:small}, the values $0$ and $2$ lie in $\fP$, $3$ is not a square in $K$, and $\varphi^{\pm2}$ are roots of $T^2-3T+1$, which has no root congruent to $1$ modulo $\fP$.
		So $z^2=1$, and hence $2a=r^2+s^2$.
		Since $\alpha=R(e)=a\,t(e)^2$, we conclude that $2\alpha=\bigl(r\,t(e)\bigr)^2+\bigl(s\,t(e)\bigr)^2$, as desired.
	\end{proof}

	\begin{remark}
		Note that in this section, we used the universality of $L$ only to represent $1$, every indecomposable that is not a square, the element $\Delta$, and, in the proof of Proposition \ref{prop:three}, the elements $3$ and $3a$.
	\end{remark}
	\section{Fields with ramified \texorpdfstring{$2$}{2}}\label{sec:dyadic}

	In this section we will show that among the fields with $U^+=U^2$ in which $2$ is ramified, only $\Q(\sqrt2)$ has the properties of the indecomposables established in Propositions \ref{prop:resP} and \ref{prop:three}.
	We will first prove that these properties force $e_\fP=2$.

	\begin{lemma}\label{lem:Dy:e}
		Let $\fP$ be a dyadic prime with $e_\fP\ge2$, and assume that every indecomposable $\alpha$ that is not a square lies in $\fP$ and that $2\alpha$ is a sum of two squares of elements of $\OO$ for every such $\alpha$.
		Then $e_\fP=2$.
	\end{lemma}

	\begin{proof}
		Suppose that $e_\fP\ge3$, and let $x$ be $S_\fP$-minimal.
		It is indecomposable by Lemma \partref{lem:minimizer}{a}, and it is not a square, as $x\notin S_\fP$.
		By our assumption, $2x=r^2+s^2$ with $r,s\in\OO$.
		Put $y=r+s$.
		As $y^2=2x+2rs\in2\OO\subset\fP^3$, we have $y\in\fP^2$.
		Hence $s\equiv-r\pmod{\fP^2}$, and so $rs\equiv-r^2\equiv r^2\pmod{\fP^2}$, as $2\in\fP^2$.
		So $rs\in S_\fP$, and the elements $x\pm rs=(r\pm s)^2/2$ do not lie in $S_\fP$.
		In particular, they are nonzero, and so they lie in $\OO^+\setminus S_\fP$.
		As $\Tr(x+rs)+\Tr(x-rs)=2\Tr(x)=2q_{S_\fP}$, both of them are $S_\fP$-minimal.
		So $x+rs=y^2/2$ is $S_\fP$-minimal with $y\in\fP^2$.
		By our assumption, Lemma \ref{lem:localtranslate} applies with $G=\fP/\fP^2$, and it gives a contradiction.
	\end{proof}

	For the remainder of the section, we will use the following hypothesis on $K$:
	\begin{itemize}
		\hypitem{(Dy)}{hyp:Dy} $U^+=U^2$, $2\alpha$ is a sum of two squares of elements of $\OO$ for every indecomposable $\alpha$, and there is a dyadic prime $\fP$ of $K$ with $e_\fP\ge2$ that contains every nonunit indecomposable.
	\end{itemize}
	Note that $\Q(\sqrt2)$ satisfies \ref{hyp:Dy} with $\fP=\sqrt2\,\OO$, as $1+\sqrt2$ has norm $-1$, the indecomposables are precisely the elements of $U^2\cup(2+\sqrt2)U^2$ (see, e.g., \cite[\S2.1]{KK26}), and $2(2+\sqrt2)=1+(1+\sqrt2)^2$.
	Under \ref{hyp:Dy}, we fix such a dyadic prime $\fP$.
	As $U^+=U^2$, every indecomposable that is not a square lies in $\fP$, and every indecomposable that does not lie in $\fP$ is the square of a unit.

	\begin{lemma}\label{lem:Dy:sqrt2}
		Assume \textup{\ref{hyp:Dy}}. Then $\sqrt2\in K$.
	\end{lemma}

	\begin{proof}
		By Lemma \ref{lem:Dy:e}, we have $v_\fP(2)=e_\fP=2$.
		Let $\gamma$ be as in Lemma \partref{lem:leastnorm}{b}.
		By \ref{hyp:Dy}, we have $2\gamma=r^2+s^2$ with $r,s\in\OO$.
		As $(r+s)^2=2\gamma+2rs\in\fP$, we have $r\equiv s\pmod\fP$.
		If $r,s\in\fP$, then $v_\fP(2\gamma)=3<4\le v_\fP(2rs)$, so $v_\fP((r+s)^2)=3$, which is impossible.
		Hence $r,s\notin\fP$, so $v_\fP(2rs)=2<3=v_\fP(2\gamma)$, and $v_\fP((r\pm s)^2)=v_\fP(2\gamma\pm2rs)=2$, i.e., $v_\fP(r\pm s)=1$.

		Put $h_\pm=(r\pm s)^2/2=\gamma\pm rs$.
		Then $h_\pm\in\OO^+$ and $v_\fP(h_\pm)=0$.
		As $U^+=U^2$, there are $\eta_\pm\in U$ with $\eta_\pm(r\pm s)\succ0$.
		
		Since $v_\fP(\eta_\pm(r\pm s))=1$, the minimality of $\N(\gamma)$ means that $|\N(r\pm s)|\ge\N(\gamma)$, i.e., $\N(h_\pm)\ge\N(\gamma)^2/2^d$.
		On the other hand, $h_+h_-=\gamma^2-r^2s^2=\bigl((r^2-s^2)/2\bigr)^2$, and $|\sigma(r^2-s^2)|<\sigma(r^2+s^2)=2\sigma(\gamma)$ for every embedding $\sigma$, as $r,s\ne0$.
		Hence $\N(h_+)\N(h_-)<\N(\gamma)^2$, and so $\N(h_\pm)<2^d$.
		By Lemma \ref{lem:basic}, $h_+$ and $h_-$ are indecomposable, and as they do not lie in $\fP$, they are squares of units by \ref{hyp:Dy}.
		So $(r+s)^2=2u^2$ for some $u\in U$, and $\sqrt2=\pm(r+s)/u\in K$.
	\end{proof}

	To prove that $K=\Q(\sqrt2)$, we will apply Lemma \ref{lem:Z} to a $\Q(\sqrt2)$-minimal element that is a square.
	The following lemma shows that such an element exists.

	\begin{lemma}\label{lem:Dy:square}
		Assume \textup{\ref{hyp:Dy}} and $K\ne\Q(\sqrt2)$.
		Then there is $y\in\OO$ such that $y^2$ is $\Q(\sqrt2)$-minimal.
	\end{lemma}

	\begin{proof}
		Put $F=\Q(\sqrt2)$, which is a subfield of $K$ by Lemma \ref{lem:Dy:sqrt2}.
		Let $x$ be $F$-minimal.
		If $x$ is a square, there is nothing to prove, so assume that it is not.
		It is indecomposable by Lemma \partref{lem:minimizer}{a}, so \ref{hyp:Dy} gives $2x=r^2+s^2$ with $r,s\in\OO$.
		Put $r'=(r+s)/\sqrt2$ and $s'=(r-s)/\sqrt2$.
		As $r'^2=x+rs$ and $s'^2=x-rs$ lie in $\OO$, we have $r',s'\in\OO$, and $r'^2+s'^2=r^2+s^2=2x$.
		If $r^2,s^2\notin F$, then $r^2,s^2\in\OO^+\setminus F$, and as $\Tr(r^2)+\Tr(s^2)=2\Tr(x)=2q_F$, the square $r^2$ is $F$-minimal.
		Similarly, if $r'^2,s'^2\notin F$, then $r'^2$ is $F$-minimal.

		Otherwise, after swapping $r$ and $s$ if necessary, which replaces $s'$ by $-s'$, and then replacing $s$ by $-s$ if necessary, which swaps $r'$ and $s'$, we may assume that $r^2,r'^2\in F$.
		As $s=\sqrt2\,r'-r$, we get $x=(r^2+s^2)/2=r^2+r'^2-\sqrt2\,rr'$.
		As $x\notin F$, this gives $rr'\notin F$, while $(rr')^2\in F$.
		Thus $rr'$ has relative trace $0$ in the quadratic extension $F(rr')/F$, and hence $\Tr(rr')=\Tr(\sqrt2\,rr')=0$.
		Therefore $y=r+r'$ satisfies $y^2=r^2+r'^2+2rr'\notin F$ and $\Tr(y^2)=\Tr(r^2+r'^2)=\Tr(x)$, i.e., $y^2$ is $F$-minimal.
	\end{proof}

	\begin{proposition}\label{prop:dyadic}
		Assume \textup{\ref{hyp:Dy}}. Then $K=\Q(\sqrt2)$.
	\end{proposition}

	\begin{proof}
		By Lemma \ref{lem:Dy:sqrt2}, $F=\Q(\sqrt2)$ is a subfield of $K$.
		Suppose that $K\ne F$, and let $y$ be as in Lemma \ref{lem:Dy:square}, so that $x=y^2$ is $F$-minimal and $y\notin F$.
		By Lemma \ref{lem:Z} with $S=F$, applied to $x=y^2$ with $A=1$ and to $x=(\sqrt2\,y)^2/2$ with $A=2$, the four elements $y^2\pm y+1$ and $y^2\pm\sqrt2\,y+1$ are totally positive and are not units, and there are $b\in\{\pm1\}$ and $c\in\{\pm\sqrt2\}$ such that each of
		\[
		z_b=y^2+by+1,\qquad z_c=y^2+cy+1
		\]
		lies in $F$ or is indecomposable.
		By \ref{hyp:Dy}, each of $z_b,z_c$ that does not lie in $F$ lies in $\fP$.
		As $z_b-z_c=(b-c)y$ with $b\ne c$ and $y\notin F$, the elements $z_b$ and $z_c$ do not both lie in $F$.
		Moreover, $b-c\in\{\pm1\pm\sqrt2\}$ is a unit, so every prime ideal containing $z_b$ and $z_c$ contains $y$, and hence also $z_b-y^2-by=1$.
		So $z_b$ and $z_c$ do not both lie in $\fP$.
		Hence one of them, say $z$, lies in $\fP$ and not in $F$, and the other one, say $z'$, lies in $F$ and not in $\fP$.

		As $z\in\fP$ and $z-1\in y\OO$, we have $y\notin\fP$.
		So $x=y^2$, which is indecomposable by Lemma \partref{lem:minimizer}{a}, does not lie in $\fP$, and it is a unit by \ref{hyp:Dy}.
		Write $z'=y^2+ky+1$ with $k\in\{b,c\}$, and put $u=z'-1=y^2+ky\in\OO\cap F=\OO_F$.
		As $y\notin F$, the minimal polynomial of $y$ over $F$ is $T^2+kT-u$, so $u=-\N_{F(y)/F}(y)$ is a unit of $\OO_F$.
		As $\OO_F/\sqrt2\,\OO_F\cong\mathbb F_2$, we get $u\equiv1\pmod{\sqrt2\,\OO_F}$, so $z'=u+1\in\sqrt2\,\OO_F\subset\fP$, a contradiction.
	\end{proof}

	\begin{corollary}\label{cor:sqclass}
		Assume that $\sqrt2\in K$ and that every indecomposable has the form $x^2$ or $(2+\sqrt2)x^2$ with $x\in\OO$.
		Then $K=\Q(\sqrt2)$.
	\end{corollary}

	\begin{proof}
		Put $\lambda=2+\sqrt2$, and let $\fP$ be a dyadic prime of $K$.
		Every $u\in U^+$ is indecomposable by Lemma \ref{lem:basic}, and it does not have the form $\lambda x^2$, as $\lambda$ is not a unit, so it is a square.
		Hence $U^+=U^2$.
		As $\sqrt2\in K$, we have $e_\fP\ge2$ and $\lambda=\sqrt2\,(1+\sqrt2)\in\fP$.
		So every indecomposable that is not a square lies in $\fP$, and by Lemma \ref{lem:notinP}, every indecomposable that does not lie in $\fP$ is a unit.
		Finally, $2x^2=x^2+x^2$ and $2\lambda x^2=x^2+\bigl((1+\sqrt2)x\bigr)^2$ for $x\in\OO$.
		So \ref{hyp:Dy} holds, and Proposition \ref{prop:dyadic} gives $K=\Q(\sqrt2)$.
	\end{proof}

	\section{Square and half-square indecomposables}\label{sec:HS}

	In this section we will determine the fields in which every indecomposable is a square or a half-square.
	We will use the following hypothesis on $K$:
	\begin{itemize}
		\hypitem{(Sq)}{hyp:Sq} $K \neq \Q$, and every indecomposable element of $\OO^+$ is a square or a half-square.
	\end{itemize}
	We will first show that under \ref{hyp:Sq} every indecomposable is a unit (Proposition \ref{prop:descent}).

	\begin{lemma}\label{lem:HSunits}
		Assume \textup{\ref{hyp:Sq}}. Then $|U^+/U^2|\le2$. If $|U^+/U^2|=2$ and $\eps\in U^+\setminus U^2$, then:
		\begin{enumerate}
			\itemlab{(a)}{lem:HSunits:a} There is $h\in\OO$ with $h^2=2\eps$, and $e_\fP$ is even for every dyadic prime $\fP$.
			\itemlab{(b)}{lem:HSunits:b} Every indecomposable element of $\OO^+$ has the form $x^2$ or $\eps x^2$ with $x\in\OO$.
		\end{enumerate}
	\end{lemma}

	\begin{proof}
		Every $u \in U^+$ is indecomposable by Lemma \ref{lem:basic}, so $u$ or $2u$ is a square by \ref{hyp:Sq}.
		If $u,u'\in U^+$ satisfy $2u=x^2$ and $2u'=x'^2$ with $x,x'\in\OO$, then $uu' = (xx'/2)^2$, so $u' \in uU^2$.
		Hence $U^+$ consists of at most $2$ classes modulo $U^2$.

		Now suppose that $|U^+/U^2|=2$, and let $\eps\in U^+\setminus U^2$.
		Then $\eps$ is not a square, so $2\eps=h^2$ for some $h\in\OO$ by the previous paragraph.
		As $(h)^2=(2)$, we have $e_\fP=2v_\fP(h)$ for every dyadic prime $\fP$, which proves \ref{lem:HSunits:a}.
		If $x \in \OO$ and $x^2/2 \in \OO$, then $x^2/2 = \eps (x/h)^2$, where $x/h \in \OO$ as $(x/h)^2=\eps^{-1}x^2/2\in\OO$.
		So every half-square is $\eps$ times a square, and \ref{lem:HSunits:b} follows from \ref{hyp:Sq}.
	\end{proof}

	\begin{proposition}\label{prop:descent}
		Assume \textup{\ref{hyp:Sq}}. Then every indecomposable is a unit.
	\end{proposition}

	\begin{proof}
	
		We will begin by showing that every indecomposable $\alpha$ can be written as $\alpha=ur^2$ with $u\in U^+$ and $r\in\OO$ such that $\eta r\succ0$ for some $\eta\in U$.
		By Lemma \ref{lem:HSunits}, we have either $U^+=U^2$ or $|U^+/U^2|=2$.

		\emph{Case $U^+=U^2$.}
		Let us show that $2$ is unramified in $K$.
		Otherwise, let $\fP$ be a dyadic prime with $e_\fP\ge2$, and let $\gamma$ be as in Lemma \partref{lem:leastnorm}{b}.
		As $v_\fP(\gamma)=1$, the element $\gamma$ is not a square, so $\gamma=y^2/2$ with $y\in\OO$, and $1=2v_\fP(y)-e_\fP$ shows that $e_\fP$ is odd, and hence $e_\fP\ge3$.
		Similarly, every indecomposable $\alpha$ that is not a square has the form $\alpha=y^2/2$ with $y\in\OO$, so $v_\fP(\alpha)=2v_\fP(y)-e_\fP$ is odd, and hence $\alpha\in\fP$, and $2\alpha=y^2+0^2$.
		So Lemma \ref{lem:Dy:e} gives $e_\fP=2$, a contradiction.

		So $e_\fP=1$ for every dyadic prime $\fP$.
		If an indecomposable $\alpha$ were a half-square, $\alpha=y^2/2$ with $y\in\OO$, then $v_\fP(\alpha)=2v_\fP(y)-1\ge0$ would give $v_\fP(y)\ge1$ and $v_\fP(\alpha)\ge1$ for every dyadic prime $\fP$.
		So $\alpha/2\in\OO^+$, and $\alpha=\alpha/2+\alpha/2$ would be decomposable.
		Hence every indecomposable is a square $r^2$, and we can take $u=1$, as there are units of all signatures.

		\emph{Case $|U^+/U^2|=2$.}
		Let $\eps$ and $h$ be as in Lemma \ref{lem:HSunits}.
		Let us first show that $e_\fP=2$ for every dyadic prime $\fP$.
		As $e_\fP$ is even by Lemma \partref{lem:HSunits}{a}, it suffices to exclude $e_\fP\ge4$.
		Suppose that $e_\fP\ge4$, and let $x$ be $S_\fP$-minimal.
		Then $x$ is indecomposable by Lemma \partref{lem:minimizer}{a}, and it is not a square, as $x\notin S_\fP$.
		So $x=\eps r^2=(hr)^2/2$ with $r\in\OO$ by Lemma \partref{lem:HSunits}{b}, and $hr\in\fP^2$, as $v_\fP(h)=e_\fP/2\ge2$.
		Also by Lemma \partref{lem:HSunits}{b}, the residue modulo $\fP^2$ of every indecomposable that is not a square lies in
		\[
		G=\{\eps t^2\bmod\fP^2:t\in\OO\},
		\]
		which is an additive subgroup of $\OO/\fP^2$, as squaring is additive on $\OO/\fP^2$.
		Moreover, $1\notin G$, since if $\eps t^2\equiv1\pmod{\fP^2}$ with $t\in\OO$, then $t\notin\fP$, so there is $t'\in\OO$ with $tt'\equiv1\pmod{\fP^2}$, and $x=\eps r^2\equiv(t'r)^2\pmod{\fP^2}$ would lie in $S_\fP$.
		So Lemma \ref{lem:localtranslate} applies, and it contradicts $x=(hr)^2/2$ with $hr\in\fP^2$.
		Hence $e_\fP=2$ and $v_\fP(h)=1$ for every dyadic prime $\fP$.

		Next, no unit has the same signature as $h$.
		Otherwise, there is $\eta\in U$ with $\eta h\succ0$, and $\eta h$ is indecomposable by Lemma \ref{lem:basic}, as $\N(\eta h)=|\N(h)|=2^{d/2}<2^d$.
		But $v_\fP(\eta h)=1$ is odd for every dyadic prime $\fP$, contrary to Lemma \partref{lem:HSunits}{b}.

		Now let $\alpha$ be an indecomposable.
		By Lemma \partref{lem:HSunits}{b}, we have $\alpha=\eps^jr^2$ with $j\in\{0,1\}$ and $r\in\OO$.
		Suppose now that no unit has the same signature as $r$, so that the signatures of $h$ and $r$ both lie outside the subgroup of unit signatures.
		As the signatures of units form a subgroup of index $2$ in $\{\pm1\}^d$ (see Section \ref{sec:prelim}), the signature of $hr$ then lies in this subgroup, so there is $\eta\in U$ with $\eta hr\succ0$.
		As $(\eta hr)^2/h^2=\eta^2\eps^{-j}\alpha$ is indecomposable and $2/h^2=\eps^{-1}\in\OO$, the element $\eta hr$ is indecomposable by Lemma \ref{lem:SU} with $A=h^2$.
		So for every dyadic prime $\fP$, the valuation $v_\fP(\eta hr)=1+v_\fP(r)$ is even by Lemma \partref{lem:HSunits}{b}, and hence $v_\fP(\alpha)=2v_\fP(r)\ge2=e_\fP$.
		Then $\alpha/2\in\OO^+$, and $\alpha$ is decomposable, a contradiction.
		So we can take $u=\eps^j$.

		\medskip

		Finally, suppose that there is a nonunit indecomposable, and let $\alpha=ur^2$ be such an indecomposable of minimal norm, with $u$, $r$, and $\eta$ as above.
		As $(\eta r)^2=\eta^2u^{-1}\alpha$ is indecomposable, $\eta r$ is indecomposable by Lemma \ref{lem:SU} with $A=1$.
		But $\eta r$ is not a unit, and $\N(\eta r)=\N(\alpha)^{1/2}<\N(\alpha)$, a contradiction.
	\end{proof}

	\begin{theorem}\label{thm:HS}
		Every indecomposable is a square or a half-square if and only if $K=\Q$, $\Q(\sqrt3)$, or $\Q(\sqrt5)$.
	\end{theorem}

	\begin{proof}
		The ``if'' part is a direct check using the description of indecomposables in real quadratic fields by Dress--Scharlau \cite[Theorem 2]{DS}.
		Indeed, over $\Q$ the only indecomposable is $1$, over $\Q(\sqrt5)$ the indecomposables are squares of units, and over $\Q(\sqrt3)$ they are the powers of the unit $2+\sqrt3$, whose odd powers are half-squares because $2(2+\sqrt3)=(1+\sqrt3)^2$.
		
		Conversely, suppose that every indecomposable is a square or a half-square and that $K\ne\Q$, i.e., that \ref{hyp:Sq} holds.
		By Proposition \ref{prop:descent}, every indecomposable is a unit.
		We will apply Lemma \ref{lem:Z} to $F$-minimal elements, first for $F=\Q$ and then for $F=\Q(\sqrt3)$ or $F=\Q(\sqrt5)$.

		For a subfield $F\subsetneq K$, let $x$ be $F$-minimal.
		Then $x$ is indecomposable by Lemma \partref{lem:minimizer}{a}, so it is a unit, and $x=y^2/A$ with $A\in\{1,2\}$ and $y\in\OO$.
		By Lemma \ref{lem:Z} with $S=F$, one of $x\pm y+1$ lies in $F$, as these elements are not units.
		Replacing $y$ by $-y$ if necessary, we may assume that $u=x+y$ lies in $F$.
		As $x\notin F$, we have $y\notin F$, and $y^2=Ax=A(u-y)$ shows that $T^2+AT-Au$ is the minimal polynomial of $y$ over $F$.
		Hence
		\begin{equation}\label{eq:normtrace}
			\N_{F(y)/F}(x)=u^2\qquad\text{and}\qquad\Tr_{F(y)/F}(x)=2u+A.
		\end{equation}
		In particular, $u$ is a unit of $\OO_F$, and $u\succ-A/4$, as the discriminant $A^2+4Au$ is totally positive.
		Moreover, for every $\eta\in U_F^+$, we have $\eta x\in\OO^+\setminus F$, so $\Tr(x)\le\Tr(\eta x)$, where $\Tr(\eta x)=\frac{[K:F]}2\Tr_{F/\Q}\bigl(\eta(2u+A)\bigr)$ by \eqref{eq:normtrace}.

		For $F=\Q$, we get $u=1$, as $u=\pm1$ and $u>-A/4\ge-1/2$.
		So $y$ is a root of $T^2+AT-A$, and hence $\sqrt5\in K$ if $A=1$, and $\sqrt3\in K$ if $A=2$.

		Next, $\sqrt3$ and $\sqrt5$ do not both lie in $K$.
		Otherwise, the element $\xi=\varphi(\sqrt3+\sqrt5)$ is totally positive, as $\varphi$ and $\sqrt3+\sqrt5$ have the sign of $\sqrt5$ at every embedding.
		As $\varphi$ is a unit, the norm of $\xi$ from $\Q(\sqrt3,\sqrt5)$ to $\Q$ is $\bigl((\sqrt5+\sqrt3)(\sqrt5-\sqrt3)\bigr)^2=4$, so $\N(\xi)=4^{d/4}<2^d$.
		So $\xi$ is a nonunit indecomposable by Lemma \ref{lem:basic}, a contradiction.

		\medskip

		Let $F$ be the field $\Q(\sqrt3)$ or $\Q(\sqrt5)$ that is contained in $K$, and suppose that $F\ne K$.
		We choose $x$, $y$, $A$, and $u$ as above for this $F$.
		If $u$ is totally negative, then $\N_{F/\Q}(u)=1$, so $u$ has a conjugate at most $-1$, contrary to $u\succ-A/4$.
		If $u$ is neither totally positive nor totally negative, then $\N_{F/\Q}(u)=-1$, so $F=\Q(\sqrt5)$, as every unit of $\Q(\sqrt3)$ has norm $1$.
		The units of $\Q(\sqrt5)$ of norm $-1$ are $\pm\varphi^n$ with odd $n\in\Z$, and the conjugate of $\varphi^n$ is $-\varphi^{-n}$.
		As $u\succ-1/2$, we get $u=\varphi^n$ or $u=-\varphi^{-n}$ with an odd $n\ge3$.
		Put $\eta=\varphi^{-2}$ in the first case and $\eta=\varphi^2$ in the second case.
		Since $\varphi^{-2}-1=-\varphi^{-1}$ and $\varphi^2-1=\varphi$, the element $(\eta-1)(2u+A)$ equals $-2\varphi^{n-1}-A\varphi^{-1}$ in the first case and $-2\varphi^{1-n}+A\varphi$ in the second case.
		As $\Tr_{F/\Q}(\varphi)=1$, $\Tr_{F/\Q}(\varphi^{-1})=-1$, and $\varphi^{1-n}$ is the conjugate of $\varphi^{n-1}$, we get in both cases
		\[
		\Tr(\eta x)-\Tr(x)=\frac{[K:F]}2\Tr_{F/\Q}\bigl((\eta-1)(2u+A)\bigr)=\frac{[K:F]}2\bigl(A-2\Tr_{F/\Q}(\varphi^{n-1})\bigr)<0,
		\]
		as $\Tr_{F/\Q}(\varphi^{n-1})\ge\Tr_{F/\Q}(\varphi^2)=3$.
		This contradicts the minimality of $\Tr(x)$.

		Thus $u$ is totally positive, and so $\N_{F/\Q}(u)=1$ and $\Tr_{F/\Q}(u^{-1})=\Tr_{F/\Q}(u)$.
		Taking $\eta=u^{-1}$, we get
		\[
		\Tr(x)-\Tr(x/u)=\frac{[K:F]}2\bigl(\Tr_{F/\Q}(2u+A)-\Tr_{F/\Q}(2+Au^{-1})\bigr)=\frac{[K:F]}2(2-A)\bigl(\Tr_{F/\Q}(u)-2\bigr).
		\]
		We will now show that we may assume that $u=1$.
		If $A=1$, the minimality of $\Tr(x)$ gives $\Tr_{F/\Q}(u)\le2$, so $\Tr(u)\le d$, and hence $u=1$ by Lemma \ref{lem:basic}.
		If $A=2$, then $x/u$ is also $F$-minimal.
		If moreover $u=\mu^2$ with $\mu\in U_F$, then Lemma \ref{lem:Z} for $x/u=(y/\mu)^2/2$ shows as above that $x/u\pm y/\mu\in F$ for one of the signs.
		As $x/u=1-y/u$ and $y\notin F$, this gives $\pm1/\mu=1/u$, so $u=\pm\mu$, and hence $u=1$.
		Otherwise $u$ is not the square of a unit, so $F=\Q(\sqrt3)$, as $U_F^+=U_F^2$ for $F=\Q(\sqrt5)$.
		Then $u$ is an odd power of $2+\sqrt3$, and $2(2+\sqrt3)=(1+\sqrt3)^2$ shows that $2u$ is a square in $F$.
		So $x/u=y^2/(2u)$ is an $F$-minimal square of an element of $\OO$, and applying the $A=1$ case to this new element gives $u=1$ for its associated parameter.

		So we may assume that $u=1$.
		Then $y$ is a root of $T^2+AT-A$, so $y\in\Q(\sqrt5)$ if $A=1$, and $y\in\Q(\sqrt3)$ if $A=2$.
		As $y\notin F$ and $K$ does not contain both $\sqrt3$ and $\sqrt5$, this is impossible.
		Hence $K=F$.
	\end{proof}

	We can now resolve the case $U^+\ne U^2$, even for lattices that need not be free.

	\begin{corollary}\label{cor:lattices}
		Assume that $U^+\ne U^2$ and that some classical ternary lattice over $\OO$ is universal. Then $K=\Q(\sqrt3)$.
	\end{corollary}

	\begin{proof}
		A universal ternary lattice can exist only if $|U^+/U^2|\le2$; see, e.g., \cite[Proposition 2.5]{KKK} (which does not assume freeness).
		So $|U^+/U^2|=2$, and \cite[Theorem 3.1(d)]{KKK} (which also does not assume freeness) shows that every indecomposable element of $\OO^+$ is a square or a half-square.
		By Theorem \ref{thm:HS}, $K=\Q$, $\Q(\sqrt3)$, or $\Q(\sqrt5)$, and of these fields only $\Q(\sqrt3)$ has $U^+ \neq U^2$.
	\end{proof}

	We are now ready to prove our main theorem.
    
	\begin{proof}[Proof of Theorem $\ref{thm:main}$]
		Each of $\Q(\sqrt2)$, $\Q(\sqrt3)$, and $\Q(\sqrt5)$ admits a universal classical ternary form by \cite[Theorem 1.1]{CKR}.
		Conversely, assume that some classical ternary form is universal over $K$.
		If $2$ is unramified in $K$, then $K=\Q(\sqrt5)$ by \cite[Theorem 1.1]{KKK}, and if $U^+\ne U^2$, then $K=\Q(\sqrt3)$ by Corollary \ref{cor:lattices}.
		So assume that $U^+=U^2$ and that $2$ is ramified in $K$, and let $\fP$ be a dyadic prime with $e_\fP\ge2$.
		Hypothesis \ref{hyp:L} holds for the free lattice $L$ corresponding to our form.
		By Proposition \partref{prop:resP}{b}, every indecomposable that does not lie in $\fP$ is a unit, and by Proposition \ref{prop:three}, $2\alpha$ is a sum of two squares of elements of $\OO$ for every indecomposable $\alpha$.
		So \ref{hyp:Dy} holds, and $K=\Q(\sqrt2)$ by Proposition \ref{prop:dyadic}.
	\end{proof}

	\printbibliography

@book{OM,
  author    = {O'Meara, O. T.},
  title     = {{Introduction to Quadratic Forms}},
  series    = {Grundlehren der mathematischen Wissenschaften},
  number    = {117},
  publisher = {Springer-Verlag},
  address   = {Berlin},
  year      = {1963},
}

@article{OM80,
  author  = {O'Meara, O. T.},
  title   = {On indecomposable quadratic forms},
  journal = {J. Reine Angew. Math.},
  volume  = {317},
  year    = {1980},
  pages   = {120--156}
}

@article{CKR,
  author       = {Chan, W. K. and Kim, M.-H. and Raghavan, S.},
  title        = {Ternary universal integral quadratic forms over real quadratic fields},
  journaltitle = {Japanese Journal of Mathematics. New Series},
  shortjournal = {Japan. J. Math.},
  volume       = {22},
  date         = {1996},
  pages        = {263--273},
}

@article{DS,
  author       = {Dress, A. and Scharlau, R.},
  title        = {Indecomposable totally positive numbers in real quadratic orders},
  journaltitle = {Journal of Number Theory},
  shortjournal = {J. Number Theory},
  volume       = {14},
  date         = {1982},
  pages        = {292--306},
}

@article{EMP,
  author       = {Edgar, H. M. and Mollin, R. A. and Peterson, B. L.},
  title        = {Class groups, totally positive units, and squares},
  journaltitle = {Proceedings of the American Mathematical Society},
  shortjournal = {Proc. Amer. Math. Soc.},
  volume       = {98},
  number       = {1},
  date         = {1986},
  pages        = {33--37},
}

@article{HKK,
  author       = {Hsia, J. S. and Kitaoka, Y. and Kneser, M.},
  title        = {Representations of positive definite quadratic forms},
  journaltitle = {Journal f{\"u}r die reine und angewandte Mathematik},
  shortjournal = {J. Reine Angew. Math.},
  volume       = {301},
  date         = {1978},
  pages        = {132--141},
}

@article{KKK,
  author       = {Kala, V. and Kramer, K. and Kr{\'a}sensk{\'y}, J.},
  title        = {{Kitaoka}'s conjecture and sums of squares},
  journaltitle = {Bulletin of the London Mathematical Society},
  shortjournal = {Bull. Lond. Math. Soc.},
  volume       = {58},
  date         = {2026},
  eid          = {e70445},
  pagetotal    = {21},
  eprint       = {2510.19545},
  eprinttype   = {arXiv},
}

@online{KK26,
  author       = {Kramer, K. and Kr{\'a}sensk{\'y}, J.},
  title        = {Non-universality of ternary quadratic forms over fields containing {$\sqrt{2}$}},
  date         = {2026-01},
  eprint       = {2601.15568},
  eprinttype   = {arXiv},
}

@online{KY24,
  author       = {Kala, V. and Yatsyna, P.},
  title        = {Even better sums of squares over quintic and cyclotomic fields},
  date         = {2024},
  eprint       = {2402.03850},
  eprinttype   = {arXiv},
}

@article{Si-tr,
  author       = {Siegel, C. L.},
  title        = {The trace of totally positive and real algebraic integers},
  journaltitle = {Annals of Mathematics. Second Series},
  shortjournal = {Ann. of Math.},
  volume       = {46},
  date         = {1945},
  pages        = {302--312},
}

@article{Si-sq,
  author       = {Siegel, C. L.},
  title        = {Sums of {$m$}-th powers of algebraic integers},
  journaltitle = {Annals of Mathematics. Second Series},
  shortjournal = {Ann. of Math.},
  volume       = {46},
  date         = {1945},
  pages        = {313--339},
}

@article{EK,
  author       = {Earnest, A. G. and Khosravani, A.},
  title        = {Universal positive quaternary quadratic lattices over totally real number fields},
  journaltitle = {Mathematika},
  volume       = {44},
  number       = {2},
  date         = {1997},
  pages        = {342--347},
  doi          = {10.1112/S0025579300012651},
}

@article{KY23,
  author       = {Kala, V. and Yatsyna, P.},
  title        = {On {Kitaoka}'s conjecture and lifting problem for universal quadratic forms},
  journaltitle = {Bulletin of the London Mathematical Society},
  shortjournal = {Bull. Lond. Math. Soc.},
  volume       = {55},
  number       = {2},
  date         = {2023},
  pages        = {854--864},
  doi          = {10.1112/blms.12762},
}

@online{Krasensky,
  author       = {Kr{\'a}sensk{\'y}, J.},
  title        = {{Kitaoka}'s conjecture for real cyclotomic fields},
  date         = {2026},
  eprint       = {2610.15109},
  eprinttype   = {arXiv},
  eprintclass  = {math.NT},
}

@article{KalaSurvey,
  author       = {Kala, V.},
  title        = {Universal quadratic forms and indecomposables in number fields: a survey},
  journaltitle = {Communications in Mathematics},
  shortjournal = {Commun. Math.},
  volume       = {31},
  number       = {2},
  date         = {2023},
  pages        = {81--114},
  doi          = {10.46298/cm.10896},
}

@article {BK,
    AUTHOR = {Blomer, Valentin and Kala, V\'it\v ezslav},
     TITLE = {Number fields without {$n$}-ary universal quadratic forms},
   JOURNAL = {Math. Proc. Cambridge Philos. Soc.},
  FJOURNAL = {Mathematical Proceedings of the Cambridge Philosophical
              Society},
    VOLUME = {159},
      YEAR = {2015},
    NUMBER = {2},
     PAGES = {239--252},
      ISSN = {0305-0041,1469-8064},
       DOI = {10.1017/S030500411500033X},
}

@article {COh,
    AUTHOR = {Chan, Wai Kiu and Oh, Byeong-Kweon},
     TITLE = {Can we recover an integral quadratic form by representing all
              its subforms?},
   JOURNAL = {Adv. Math.},
  FJOURNAL = {Advances in Mathematics},
    VOLUME = {433},
      YEAR = {2023},
     PAGES = {Paper No. 109317},
 PAGETOTAL = {20},
       DOI = {10.1016/j.aim.2023.109317},
}

@article {HHX,
    AUTHOR = {He, Zilong and Hu, Yong and Xu, Fei},
     TITLE = {On indefinite {$k$}-universal integral quadratic forms over
              number fields},
   JOURNAL = {Math. Z.},
  FJOURNAL = {Mathematische Zeitschrift},
    VOLUME = {304},
      YEAR = {2023},
    NUMBER = {1},
     PAGES = {Paper No. 20},
 PAGETOTAL = {26},
       DOI = {10.1007/s00209-023-03280-z},
}

@article {KT,
    AUTHOR = {Kala, V\'it\v ezslav and Tinkov\'a, Magdal\'ena},
     TITLE = {Universal quadratic forms, small norms, and traces in families
              of number fields},
   JOURNAL = {Int. Math. Res. Not. IMRN},
  FJOURNAL = {International Mathematics Research Notices. IMRN},
    VOLUME = {2023},
      YEAR = {2023},
    NUMBER = {9},
     PAGES = {7541--7577},
       DOI = {10.1093/imrn/rnac073},
}

@article {XZ,
    AUTHOR = {Xu, Fei and Zhang, Yang},
     TITLE = {On indefinite and potentially universal quadratic forms over
              number fields},
   JOURNAL = {Trans. Amer. Math. Soc.},
  FJOURNAL = {Transactions of the American Mathematical Society},
    VOLUME = {375},
      YEAR = {2022},
    NUMBER = {4},
     PAGES = {2459--2480},
       DOI = {10.1090/tran/8601},
}

@article {Yat19,
    AUTHOR = {Yatsyna, Pavlo},
     TITLE = {A lower bound for the rank of a universal quadratic form with
              integer coefficients in a totally real number field},
   JOURNAL = {Comment. Math. Helv.},
  FJOURNAL = {Commentarii Mathematici Helvetici. A Journal of the Swiss
              Mathematical Society},
    VOLUME = {94},
      YEAR = {2019},
    NUMBER = {2},
     PAGES = {221--239},
      ISSN = {0010-2571,1420-8946},
   MRCLASS = {11E12 (11H06 11R80)},
  MRNUMBER = {3942784},
MRREVIEWER = {Lenny\ Fukshansky},
       DOI = {10.4171/CMH/459},
       URL = {https://doi.org/10.4171/CMH/459},
}

@article {KKP,
    AUTHOR = {Kim, Byeong Moon and Kim, Myung-Hwan and Park, Dayoon},
     TITLE = {Real quadratic fields admitting universal lattices of rank 7},
   JOURNAL = {J. Number Theory},
  FJOURNAL = {Journal of Number Theory},
    VOLUME = {233},
      YEAR = {2022},
     PAGES = {456--466},
      ISSN = {0022-314X,1096-1658},
   MRCLASS = {11E12},
  MRNUMBER = {4356860},
MRREVIEWER = {Jingbo\ Liu},
       DOI = {10.1016/j.jnt.2021.06.018},
       URL = {https://doi.org/10.1016/j.jnt.2021.06.018},
}

@article {Kim,
    AUTHOR = {Kim, Byeong Moon},
     TITLE = {Universal octonary diagonal forms over some real quadratic
              fields},
   JOURNAL = {Comment. Math. Helv.},
  FJOURNAL = {Commentarii Mathematici Helvetici},
    VOLUME = {75},
      YEAR = {2000},
    NUMBER = {3},
     PAGES = {410--414},
      ISSN = {0010-2571,1420-8946},
   MRCLASS = {11E12 (11E20)},
  MRNUMBER = {1793795},
MRREVIEWER = {Andrzej\ S\l adek},
       DOI = {10.1007/s000140050133},
       URL = {https://doi.org/10.1007/s000140050133},
}

@article{Ma,
  author  = {Maa{\ss}, H.},
  title   = {{\"U}ber die {D}arstellung total positiver {Z}ahlen des {K}{\"o}rpers {$R(\sqrt5)$} als {S}umme von drei {Q}uadraten},
  journal = {Abh. Math. Sem. Univ. Hamburg},
  volume  = {14},
  year    = {1941},
  pages   = {185--191},
}

@incollection{Bh,
  author    = {Bhargava, M.},
  title     = {On the {C}onway--{S}chneeberger fifteen theorem},
  booktitle = {Quadratic Forms and Their Applications ({D}ublin, 1999)},
  series    = {Contemp. Math.},
  number    = {272},
  publisher = {Amer. Math. Soc.},
  address   = {Providence, RI},
  year      = {2000},
  pages     = {27--37},
}

@article{KTZ,
  author  = {Kr{\'a}sensk{\'y}, J. and Tinkov{\'a}, M. and Zemkov{\'a}, K.},
  title   = {There are no universal ternary quadratic forms over biquadratic fields},
  journal = {Proc. Edinb. Math. Soc.},
  volume  = {63},
  year    = {2020},
  pages   = {861--912},
}

@misc{KKPYZ,
  author = {Kala, V. and Kr{\'a}sensk{\'y}, J. and Park, D. and Yatsyna, P. and {\.Z}mija, B.},
  title  = {{Kitaoka}'s conjecture for quadratic fields},
  note   = {arXiv:2501.19371},
  year   = {2025},
}

@article{KYZ,
  author  = {Kala, V. and Yatsyna, P. and {\.Z}mija, B.},
  title   = {Real quadratic fields with a universal quadratic form of given rank have density zero},
  journal = {Amer. J. Math.},
  note    = {to appear, arXiv:2302.12080},
}

@article{BK18,
  author  = {Blomer, V. and Kala, V.},
  title   = {On the rank of universal quadratic forms over real quadratic fields},
  journal = {Doc. Math.},
  volume  = {23},
  year    = {2018},
  pages   = {15--34},
}

@article{Ka16,
  author  = {Kala, V.},
  title   = {Universal quadratic forms and elements of small norm in real quadratic fields},
  journal = {Bull. Aust. Math. Soc.},
  volume  = {94},
  year    = {2016},
  pages   = {7--14},
}

@article{Ka23,
  author  = {Kala, V.},
  title   = {Number fields without universal quadratic forms of small rank exist in most degrees},
  journal = {Math. Proc. Cambridge Philos. Soc.},
  volume  = {174},
  year    = {2023},
  pages   = {225--231},
}

@article{Man,
  author  = {Man, S. H.},
  title   = {Minimal rank of universal lattices and number of indecomposable elements in real multiquadratic fields},
  journal = {Adv. Math.},
  volume  = {447},
  year    = {2024},
  pages   = {Paper No. 109694},
  pagetotal = {38},
}

@incollection{Kim04,
  author    = {Kim, M.-H.},
  title     = {Recent developments on universal forms},
  booktitle = {Algebraic and arithmetic theory of quadratic forms},
  series    = {Contemp. Math.},
  number    = {344},
  publisher = {Amer. Math. Soc.},
  address   = {Providence, RI},
  year      = {2004},
  pages     = {215--228},
}

@online{Hotta-I,
  author = {Hotta, Kika Zen-U L.},
  title  = {{Kitaoka}'s conjecture and universal ternary lattices, {Part I}: A degree bound and the case of {$\mathbb{Q}(\sqrt{10})$}},
  date   = {2026-11},
  note   = {DOI:10.13140/RG.2.2.21318.02883},
}

@online{Hotta-II,
  author = {Hotta, Kika Zen-U L.},
  title  = {{Kitaoka}'s conjecture and universal ternary lattices, {Part II}: The optimal degree bound for classical ternary forms},
  date   = {2026-11},
  note   = {DOI:10.13140/RG.2.2.31928.30721},
}

\end{document}